\documentclass[12pt]{article}      

\title{Golden Binomials and Carlitz Characteristic Polynomials}

\author{ Oktay K Pashaev and Merve \"Ozvatan\\Department of Mathematics\\ Izmir Institute of Technology \\ Urla-Izmir, 35430, Turkey}

\begin{document}
\newcommand{\be}{\begin{equation}}
\newcommand{\ee}{\end{equation}}
\newcommand{\bea}{\begin{eqnarray}}
\newcommand{\eea}{\end{eqnarray}}
\newcommand{\disp}{\displaystyle}
\newcommand{\la}{\langle}
\newcommand{\ra}{\rangle}

\newtheorem{thm}{Theorem}[subsection]
\newtheorem{cor}[thm]{Corollary}
\newtheorem{lem}[thm]{Lemma}
\newtheorem{prop}[thm]{Proposition}
\newtheorem{definition}[thm]{Definition}
\newtheorem{rem}[thm]{Remark}
\newtheorem{prf}[thm]{Proof}

\maketitle


\begin{abstract}
The golden binomials, introduced  in the golden quantum calculus, have expansion determined by Fibonomial coefficients and the set of simple zeros given by powers of Golden ratio. We show that these golden binomials are equivalent to Carlitz characteristic polynomials of certain matrices of binomial coefficients.  It is shown that trace invariants for powers of these matrices are determined by Fibonacci divisors, quantum calculus of which was developed very recently.
\end{abstract}

\section{Introduction}
The golden quantum calculus,  based on the Binet formula for
Fibonacci numbers $F_n$ as $q$-numbers, 
was introduced in \cite{golden}. In this calculus, the
 finite-difference $q$-derivative operator is determined by two Golden ratio bases $\varphi$ and $\varphi'$, while the golden binomial expansion, by Fibonomial coefficients. The coefficients  are expressed in terms of  Fibonacci numbers, while zeros of these binomials are given by powers of Golden ratio $\varphi$ and $\varphi'$.
It was observed that similar polynomials were introduced by Carlitz in 1965 from different reason,  as characteristic polynomials of certain matrices of binomial coefficients \cite{Carlitz}. The goal of the present paper is to show equivalence of Carlitz characterisitc polynomials with golden binomials. In addition, the proof and interpretation of main formulas for trace of powers of the matrix $A_{n+1}$ in terms of Fibonacci divisors and corresponding quantum calculus, developed recently  in \cite{FNK} would be given.

\section{Golden Binomials}

\subsection{Fibonomials and Golden Pascal Triangle}

The  binomial coefficients  defined by
\be { n \brack k}_F= \frac{[n]_{F}!}{[n-k]_{F}! [k]_{F}!}= \frac{F_n!}{F_{n-k}! F_k!}, \label{goldenbinom}\ee
with $n$ and $k$ being non-negative integers, $n\geq k$, are called the Fibonomials.
Using the addition formula for Fibonacci  numbers \cite{golden}, 
\be
F_{n+m} = \varphi^n F_m + {\varphi'}^m F_n
\ee
we have following expression
\be 
F_n=F_{n-k+k}=\left(-\frac{1}{\varphi}\right)^k F_{n-k}+\varphi^{n-k} F_k.
\ee
By using
\be
\varphi^n = \varphi F_n + F_{n-1}, \,\,\,\,{\varphi'}^n = \varphi' F_n + F_{n-1},\label{phin}
\ee
 it can be rewritten as follows
\bea F_n &=& F_{n-k-1} F_k + F_{n-k} F_{k+1} \nonumber \\
&=& F_{n-k} F_{k-1}+ F_{n-k+1} F_k. \eea
With the above definition (\ref{goldenbinom}) it gives recursion formula for Fibonomials in two forms,
\bea { n \brack k}_{F}&=& \frac{(-\frac{1}{\varphi})^k [n-1]_{F}!}{[k]_{F}! [n-k-1]_{F}!} + \frac{\varphi^{n-k} [n-1]_{F}!}{[n-k]_{F}![k-1]_{F}!} \nonumber \\
&=& \left(-\frac{1}{\varphi}\right)^k { n-1 \brack k}_{F} + \varphi^{n-k} { n-1 \brack k-1}_{F} \label{goldenpascal1}\\
&=& \varphi^k { n-1 \brack k}_{F} + \left(-\frac{1}{\varphi}\right)^{n-k} { n-1 \brack k-1}_{F} .\label{goldenpascal2} \eea
These formulas, for $1\leq k\leq n-1$, determine the  Golden Pascal triangle for Fibonomials  \cite{golden}.

\subsection{Golden Binomial}
The Golden Binomial is  defined  as \cite{golden},
\be (x+y)_F^n = (x+\varphi^{n-1} y)(x+\varphi^{n-2} \varphi' y)...(x+ \varphi{\varphi'}^{n-2} y)
(x+ {\varphi'}^{n-1} y)\ee
or due to $\varphi \varphi' = -1$ it is
\be (x+y)_F^n = (x+\varphi^{n-1} y)(x-\varphi^{n-3} y)...(x+ (-1)^{n-1}\varphi^{-n+1} y).\ee 
It has n-zeros at  powers of the Golden ratio 
$$\frac{x}{y}=-\varphi^{n-1},\,\,\,\, \frac{x}{y}=-\varphi^{n-3},\,\,\,\,...,\frac{x}{y}=-\varphi^{-n+1}.$$
For Golden binomial the following expansion in terms of Fibonomials is valid \cite{golden}
\bea (x+y)_F^n &=&\sum^{n}_{k=0}{ n \brack k}_{F} (-1)^{\frac{k(k-1)}{2}} x^{n-k} y^k  \nonumber \\
&=& \sum^{n}_{k=0} \frac{F_n!}{F_{n-k}! F_k!}(-1)^{\frac{k(k-1)}{2}}  x^{n-k} y^k. \label{goldenbinomexpansion}\eea
The proof is easy by induction and using recursion formulas (\ref{goldenpascal1}), (\ref{goldenpascal2}) . In terms of Golden binomials we introduce the Golden polynomials
\be P_n (x) = \frac{(x-a)_F^n}{F_n!},\ee
where $n=1,2,...$, and $P_0(x) =1$. These polynomials satisfy relations
\be D_F^x P_n(x) = P_{n-1}(x), \ee
 where the Golden derivative is defined as
 \be D_F^x P_n(x) = \frac{P_n (\varphi x) - P_n (\varphi' x) }{(\varphi - \varphi') x}.\ee
For even and odd polynomials we have different products
\be P_{2n} (x) = \frac{1}{F_{2n}!} \prod^n_{k=1} (x- (-1)^{n+k}\varphi^{2k-1} a) (x + (-1)^{n+k}\varphi^{-2k +1} a) ,\ee
\be P_{2n+1} (x) = \frac{(x - (-1)^n a)}{F_{2n+1}!} \prod^n_{k=1} (x- (-1)^{n+k}\varphi^{2k} a) (x - (-1)^{n+k}\varphi^{-2k} a) .\ee
By using (\ref{phin}) it is easy to find
\be \varphi^{2k} + \frac{1}{\varphi^{2k}} = F_{2k} + 2 F_{2k-1} ,\ee
\be \varphi^{2k+1} - \frac{1}{\varphi^{2k+1}} = F_{2k+1} + 2 F_{2k} .\ee
Then we can rewrite our polynomials in terms of Fibonacci numbers
\be P_{2n} (x) = \frac{1}{F_{2n}!} \prod^n_{k=1} (x^2 - (-1)^{n+k} (F_{2k-1} + 2 F_{2k-2})x a - a^2) ,\ee
\be P_{2n+1} (x) = \frac{(x - (-1)^n a)}{F_{2n+1}!} \prod^n_{k=1} (x^2- (-1)^{n+k}(F_{2k} + 2 F_{2k-1})x a + a^2) .\ee
The first few odd polynomials are
\be P_1(x) = (x-a),\ee
\be P_3 (x) = \frac{1}{2} (x+a)(x^2 - 3 x a + a^2 ),\ee
\be P_5 (x) = \frac{1}{2\cdot 3 \cdot 5} (x-a)(x^2 + 3 x a + a^2 )(x^2 - 7 x a + a^2 ),\ee
\be P_7 (x) = \frac{1}{2\cdot 3 \cdot 5 \cdot 8 \cdot 13} (x+a)(x^2 - 3 x a + a^2 )(x^2 + 7 x a + a^2 )(x^2 - 18 x a + a^2 ),\ee
and the even ones
\be P_2 (x) =  (x^2 -  x a - a^2 ),\ee
\be P_4 (x) = \frac{1}{2\cdot 3} (x^2 +  x a - a^2 )(x^2 - 4 x a - a^2),\ee
\be P_6 (x) = \frac{1}{2\cdot 3\cdot 5 \cdot 8} (x^2 -  x a - a^2 )(x^2 + 4 x a - a^2)(x^2 - 11 x a - a^2).\ee

\subsection{Golden analytic function}
By golden binomials in complex domain, the golden analytic function can be derived, which is complex valued  function of complex argument, not analytic in usual sense \cite{Eskisehir}. 
The complex golden binomial is defined as
\bea (x+iy)^n_F &=& (x + i\varphi^{n-1}y)(x - i\varphi^{n-3} y)... (x + i(-1)^{n-1}\varphi^{1-n}y)  \\ &=&\sum^n_{k=0}\left[\begin{array}{c}n \\ k \end{array}\right]_{F} (-1)^{\frac{k(k-1)}{2}}x^{n-k}i^k y^k.\eea
It can be generated by the golden translation
$$ E^{iy D^x_F}_F x^n = (x+iy)^n_F,$$
where
$$ E^x_F = \sum^\infty_{n=0} (-1)^{\frac{n(n-1)}{2}} \frac{x^n}{F_n!}.$$
The binomials determine the golden analytic function
$$ f(z, F) = E^{iy D^x_F}_F f(x) = \sum^\infty_{n=0} a_n \frac{(x+iy)^n_F}{F_n!},$$
satisfying the golden $\bar\partial_F$ equation
\be  \frac{1}{2}(D^x_{F} + i D^y_{-F}) f(z;F) = 0,\ee
where
$D^x_{-F} = (-1)^{x\frac{d}{dx}} D^x_F$.
For $u(x,y) = Cos_{F} (y D^x_{F}) f(x)$  and $v(x,y) = Sin_{F} (y D^x_{F}) f(x)$,
the golden Cauchy-Riemann equations are
\be D^x_{F} u(x,y) = D^y_{-F} v(x,y),\,\,\,\,D^y_{-F} u(x,y) = -D^x_{F} v(x,y),\ee
and the golden-Laplace equation is
\be (D^x_{F})^2 u(x,y) + (D^y_{-F})^2 u(x,y) = 0. \ee

\subsection{Particular Case}

The golden binomial $(x -a)^n_F$ can be also generated by the golden translation
\be
E^{-a D^x_F}_F x^n = (x-a)^n_F.
\ee
In particular case $a = 1$ we have

\be
(x-1)^m_F = (x - \varphi^{m-1}) (x + \varphi^{m-3})...  (x - (-1)^{m-1}\varphi^{-m+1}) .\label{x1}
\ee
First few binomials are
\bea
(x-1)^1_F &=& x -1 ,\\
(x-1)^2_F& =& (x -\varphi) (x - \varphi'), \\
(x-1)^3_F& =& (x -\varphi^2)(x+1) (x - {\varphi'}^2), \\
(x-1)^4_F& =& (x -\varphi^3)  (x +\varphi) (x + \varphi')     (x - {\varphi'}^3) ,
\eea
and corresponding zeros
\bea
m = 1 &\Rightarrow & x =1 \\
m = 2 & \Rightarrow& x =\varphi, x = \varphi'\\
m=3 & \Rightarrow & x =\varphi^2, x = -1, x = {\varphi'}^2 \\
m=4& \Rightarrow& x =\varphi^3,   x = -\varphi, x = -\varphi',     x = {\varphi'}^3.
\eea
For arbitrary even and odd  $n$ we have following zeros of Golden binomials
\bea
n = 2k & \Rightarrow    & (x-1)^{2k}_F : \varphi^{n-1}, {\varphi'}^{n-1}, -\varphi^{n-3}, -{\varphi'}^{n-3}, ..., \pm\varphi, \pm {\varphi'};\label{even}\\
 n = 2k+1 & \Rightarrow    & (x-1)^{2k+1}_F : \varphi^{n-1}, {\varphi'}^{n-1}, -\varphi^{n-3}, -{\varphi'}^{n-3}, ..., \pm 1.\label{odd}
\eea

\section{Carlitz Polynomials}

In Section 2 we have introduced the Golden binomials. Now we are going to relate these binomials with characteristic equations for some matrices, constructed from binomial coefficients by Carlitz \cite{Carlitz}.
\begin{definition}
We define an $(n+1) \times (n+1)$ matrix $A_{n+1}$ with binomial coefficients,
\bea
A_{n+1}=\left[{r \choose n-s} \right],
\eea
where $r,s=0,1,2,...,n .$ Here, \begin{eqnarray}
{n \choose k} =\left\{          \begin{array}{ll}
                                        \frac{n!}{(n-k)! \phantom{.} k!}, & \hbox{if k $\leq$ n;} \\
                                         0, & \hbox{k $>$ n.}
                                       \end{array}
                                     \right.
\end{eqnarray}
\end{definition}
First few matrices are,
\bea
&&{n=0} \phantom{a} \Rightarrow r=s=0 \Rightarrow \phantom{a} A_{1}=\left[{0 \choose 0} \right]=(1) \nonumber \\
&&{n=1} \phantom{a} \Rightarrow r,s=0,1 \Rightarrow \phantom{a} A_{2}=\left[{r \choose 1-s} \right]=\left(
                                       \begin{array}{cc}
                                          {0 \choose 1} & {0 \choose 0} \\
                                          {1 \choose 1} & {1 \choose 0} \\
                                       \end{array}
                                     \right)=\left(
                                       \begin{array}{cc}
                                         0 & 1 \\
                                         1 & 1 \\
                                       \end{array}
                                     \right)\nonumber \\
&&{n=2} \Rightarrow r,s=0,1,2 \Rightarrow A_{3}=\left[{r \choose 2-s} \right]=\left(
                                       \begin{array}{ccc}
                                         {0 \choose 2} & {0 \choose 1} & {0 \choose 0} \\
                                         {1 \choose 2} & {1 \choose 1} & {1 \choose 0} \\
                                         {2 \choose 2} & {2 \choose 1} & {2 \choose 0} \\
                                       \end{array}
                                     \right)=\left(
                                       \begin{array}{ccc}
                                         0 & 0 & 1 \\
                                         0 & 1 & 1 \\
                                         1 & 2 & 1 \\
                                       \end{array}
                                     \right)
 \nonumber
\eea
Continuing, the general matrix $A_{n+1}$ of order $(n+1)$ can be written as,
\begin{eqnarray}
A_{n+1}=\left(
  \begin{array}{ccccc}
   \ldots \phantom{.} 0 & 0 & 0 & 0 & 1 \\
   \ldots \phantom{.} 0 & 0 & 0 & 1 & 1 \\
   \ldots \phantom{.} 0 & 0 & 1 & 2 & 1 \\
   \ldots \phantom{.} 0 & 1 & 3 & 3 & 1 \\
   \ldots \phantom{.} 1 & 4 & 6 & 4 & 1 \\
   \phantom{.....} \vdots & \vdots & \vdots & \vdots & \vdots \\
  \end{array}
\right)_{(n+1)\times(n+1)}, \nonumber
\end{eqnarray}
where the lower triangular matrix is build from Pascal's triangle. 
We notice that trace of first few matrices $A_{n+1}$ gives Fibonacci numbers. As would be shown, it is valid for any n (Theorem $(\ref{invarianttheorem})$ equation $(\ref{invarianttheoremequation1})$) .

\begin{definition}
Characteristic polynomial of matrix $A_{n+1}$ is determined by,
\bea
Q_{n+1}(x)=\mbox{det}(x I-A_{n+1}). \label{characteristicequation}
\eea
\end{definition}
First few polynomials explicitly are

\bea
&&{n=0:}\phantom{abc} Q_{1}(x)=x - 1 ,\nonumber \\
&&{n=1:}\phantom{abc} Q_{2}(x)=\mbox{det}(x I-A_{2})=\left|
                                                                        \begin{array}{cc}
                                                                          x & -1 \\
                                                                          -1 & x-1 \\
                                                                        \end{array}
                                                                      \right|
=x^2-x-1, \nonumber \\
&&{n=2:}\phantom{abc} Q_{3}(x)=\mbox{det}(x I-A_{3})=\left|
                                                                              \begin{array}{ccc}
                                                                                x & 0 &  -1 \\
                                                                                0 & x-1 & -1 \\
                                                                                -1 & -2 & x-1 \\
                                                                              \end{array}
 \right|=x^3-2 x^2-2 x+1 ,\nonumber \\
&&{n=3:}\phantom{abc} \nonumber
\eea
\bea
Q_{4}(x)=\mbox{det}(x I-A_{4})&=&\left|
                                               \begin{array}{cccc}
                                                 x & 0 & 0 & -1 \\
                                                 0 & x & -1 & -1 \\
                                                 0 & -1 & x-2 & -1 \\
                                                -1 & -3 & -3 & x-1 \\
                                               \end{array}
                                             \right| \nonumber \\
&=&-x^4+3 x^3+6 x^2-3 x-1 . \nonumber
\eea

Corresponding eigenvalues are represented by powers of $\varphi$ and $\varphi'$;

${n=0}$\phantom{abc} $\Rightarrow$ \phantom{a} $x_1=1$,

${n=1}$\phantom{abc}$\Rightarrow$ \phantom{a} $x_1=\varphi,\quad x_2=\varphi'$,

${n=2}$\phantom{abc}$\Rightarrow$ \phantom{a} $x_1=\varphi^2,\quad x_2=-1,\quad x_3={\varphi'}^2$,

${n=3}$\phantom{abc}$\Rightarrow$ \phantom{a} $x_1=\varphi^3, x_2=-\varphi, x_3=-\varphi',x_4=\varphi'^3$.

Comparing zeros of first few characteristic polynomials, with zeros of Golden Binomial $(\ref{x1})$, we notice that they coincide. According to this, we have following conjecture.

\textbf{Conjecture:} The characteristic equation $(\ref{characteristicequation})$ of matrix $A_{n+1}$ coincides with Golden Binomial;
\bea
Q_{n+1}(x)=\mbox{det}(x I-A_{n+1})=(x-1)^{n+1}_{F}.
\eea
 To prove this conjecture, firstly we represent Golden binomials in the product form.

\begin{prop}
The Golden binomial can be written as a product,
\begin{eqnarray}
(x-1)^{n+1}_{F}=\prod_{j=0}^{n} \left(x-\varphi^j \varphi'^{n-j}\right).
\end{eqnarray}
\end{prop}
\begin{prf}
Starting from Golden binomial in product representation
\bea
(x+y)^{n}_{F} \equiv \prod_{j=0}^{n-1} \left(x-(-1)^{j-1}\phantom{.}  \varphi^{n-1}\phantom{.}  \varphi^{-2j} y   \right)
\eea
by using
\bea
\varphi^{-2j}=\left(\frac{1}{\varphi}  \right)^{2j}=\left(-\frac{1}{\varphi}  \right)^{2j}=\varphi'^{2j},
\eea
after substitution $y=-1$ we have
\bea
(x-1)^{n}_{F} \equiv \prod_{j=0}^{n-1} \left(x-(-1)^{j} \phantom{.} \varphi^{n-1}\phantom{.} \varphi'^{2j}   \right) .\nonumber
\eea
By shifting $n \rightarrow n+1$, 
\bea
(x-1)^{n+1}_{F}&=&\prod_{j=0}^{n} \left(x-(-1)^{j} \phantom{.} \varphi^{n}\phantom{.} \varphi'^{2j} \right) \nonumber \\
&=&\prod_{j=0}^{n} \left(x-(-1)^{j} \phantom{.} \varphi^{n}\phantom{.} \frac{(-1)^{2j}}{\varphi^{j} \varphi^{j}}   \right) \nonumber \\
&=&\prod_{j=0}^{n} \left(x- \varphi^{n} \left(-\frac{1}{\varphi}\right)^{j} \phantom{.} \frac{1}{\varphi^{j}}   \right) \nonumber \\
&=&\prod_{j=0}^{n} \left(x- \varphi^{n-j} \varphi'^{j} \phantom{.}    \right) \nonumber
\eea
and substituting  $j=n-m$ we get,
\bea
(x-1)^{n+1}_{F}=\prod_{m=0}^{n} \left(x- \varphi^{m} \varphi'\phantom{.}^{n-m} \phantom{.}    \right). \nonumber
\eea
The formula  shows explicitly that zeros of Golden binomial in $(\ref{even})$ and $(\ref{odd})$ are given  by powers of $\varphi$ and $\varphi'$.
\end{prf}

\begin{cor}
The eigenvalues of matrix $A_{n+1}$ are the numbers,
\bea
\varphi^n, \varphi^{n-1}\varphi', \varphi^{n-2}\varphi'^{2}, \ldots ,\varphi \phantom{.}\varphi'^{n-1}, \varphi'^{n} \label{eigenvaluesofmatrixAn+1}.
\eea
\end{cor}

As it was shown by Carlitz \cite{Carlitz},  this product formula is just characteristic equation $(\ref{characteristicequation})$ for matrix $A_{n+1}$. Since zeros of two polynomials $\mbox{det}(x I -A_{n+1})$ and $(x-1)^{n+1}_{F}$ coincide, then the conjecture is correct and we have following theorem.
\begin{thm}
Characteristic equation for combinatorial matrix $A_{n+1}$ is given by Golden binomial:
\bea
Q_{n+1}(x)=\mbox{det}(x I-A_{n+1})=(x-1)^{n+1}_{F}.
\eea
\end{thm}

\section{Powers of $A_{n+1}$ and Fibonacci Divisors}

\begin{prop}
Arbitrary $n^{th}$ power of $A_{2}$ matrix is written in terms of Fibonacci numbers,
\bea
A^{n}_{2}=\left( \begin{array}{cc}
   F_{n-1} & F_{n} \\
   F_{n} & F_{n+1} \\
  \end{array}
     \right).
\eea
\end{prop}
\begin{prf}
Proof will be done by induction. For $n=1$,
\bea
A_{2}=\left(
\begin{array}{cc}
                           0 & 1 \\
                           1 & 1
                        \end{array} \right) =\left( \begin{array}{cc}
                           F_{0} & F_{1} \\
                           F_{1} & F_{2}
                        \end{array}\right), \nonumber
\eea
and for $n=2$,
\bea
A^2_{2}=\left(
\begin{array}{cc}
                           1 & 1 \\
                           1 & 2
                        \end{array} \right) =\left( \begin{array}{cc}
                           F_{1} & F_{2} \\
                           F_{2} & F_{3}
                        \end{array}\right). \nonumber
\eea
Suppose for $n=k$,
\bea
A^{k}_{2}=\left( \begin{array}{cc}
   F_{k-1} & F_{k} \\
   F_{k} & F_{k+1} \\
  \end{array}
     \right) , \nonumber
\eea
then 
\bea
A^{k+1}_{2}&=&A^{k}_{2}\phantom{.} A_{2}=\left( \begin{array}{cc}
   F_{k-1} & F_{k} \\
   F_{k} & F_{k+1} \\
  \end{array}
     \right)\phantom{.} \left( \begin{array}{cc}
   0 & 1 \\
   1 & 1 \\
  \end{array}
     \right) \\
&=& \left( \begin{array}{cc}
   F_{k} & F_{k}+F_{k-1} \\
   F_{k+1} & F_{k}+F_{k+1} \\
  \end{array}
     \right)=\left( \begin{array}{cc}
   F_{k} & F_{k+1} \\
   F_{k+1} & F_{k+2} \\
  \end{array}
     \right).\nonumber
\eea
This result can be understood from observation that eigenvalues of matrix $A_2$ are $\varphi$ and $\varphi'$, and eigenvalues of $A^{n}_2$ are powers $\varphi^n$, $\varphi'^n$ related with Fibonacci numbers.

\end{prf}

As we have seen, eigenvalues of matrix $A_3$ are $\varphi^2, \varphi'^2, -1$. It implies that for $A^{n}_3$, eigenvalues are $\varphi^{2n}, \varphi'^{2n}, (-1)^n$, and the matrix can be expressed by Fibonacci divisor $F^{(2)}_n$ conjugate to $F_2$, due to  \cite{FNK},
\bea
(\varphi^k)^n = \varphi^k F^{(k)}_n + (-1)^{k+1}  F^{(k)}_{n-1},\label{phikn1}\\
({\varphi'}^k)^n = {\varphi'}^k F^{(k)}_n + (-1)^{k+1}  F^{(k)}_{n-1},\label{phikn2}
\eea
where $F^{(k)}_n = F_{nk}/F_k$.
\begin{prop} \label{A3powern}
Arbitrary $n^{th}$ power of $A_{3}$ matrix can be expressed in terms of Fibonacci divisors $F_{n}^{(2)}$,
\begin{eqnarray}
A^{n}_3= \frac{1}{5}\left(
  \begin{array}{ccc}
    (2F_{n}^{(2)}-3F_{n-1}^{(2)}+2(-1)^n) & (2F_{n}^{(2)}+2F_{n-1}^{(2)}+2(-1)^n) & (3F_{n}^{(2)}-2F_{n-1}^{(2)}-2(-1)^n) \\
    (F_{n}^{(2)}+F_{n-1}^{(2)}+(-1)^n) & (6F_{n}^{(2)}-4F_{n-1}^{(2)}+(-1)^n) & (4F_{n}^{(2)}-F_{n-1}^{(2)}-(-1)^n) \\
    (3F_{n}^{(2)}-2F_{n-1}^{(2)}-2(-1)^n) & (8F_{n}^{(2)}-2F_{n-1}^{(2)}-2(-1)^n) & (7F_{n}^{(2)}-3F_{n-1}^{(2)}+2(-1)^n) \\
  \end{array}
\right) \nonumber
\end{eqnarray}
\end{prop}
\begin{prf}
Let's diagonalize the matrix $A_{3}$,
\bea
\phi_{3}=\sigma^{-1}_{3}\phantom{a} A_{3} \phantom{a} \sigma_{3}, \nonumber
\eea
where $\phi_{3}$ is the diagonal matrix and
\bea
A_{3}=\sigma_{3} \phantom{.} \phi_{3} \phantom{.} \sigma^{-1}_{3}. \nonumber
\eea
Taking the $n^{th}$ power of both sides gives,
\bea
A^{n}_{3}=(\sigma_{3} \phantom{.} \phi_{3} \phantom{.} \underbrace{\sigma^{-1}_{3})\phantom{a}(\sigma_{3}}_{I} \phantom{.} \phi_{3} \phantom{.} \sigma^{-1}_{3})\phantom{a}...\phantom{a}(\sigma_{3} \phantom{.} \phi_{3} \phantom{.} \underbrace{ \sigma^{-1}_{3})\phantom{a}(\sigma_{3}}_{I} \phantom{.} \phi_{3}  \phantom{.} \sigma^{-1}_{3}) \nonumber
\eea
Therefore, 
\bea
{A^{n}_{3}=\sigma_{3} \phantom{.} \phi^{n}_{3} \phantom{.} \sigma^{-1}_{3}}. \label{A3powerintermsofphiandsigma}
\eea
By using the diagonalization principle, $\sigma_{3}$ and $\sigma^{-1}_{3}$ matrices can be obtained as,
\bea
\sigma_{3}=\frac{1}{2}\left(
\begin{array}{ccc}
  -\varphi' & \frac{4}{3} & -\varphi \\
  1 & \frac{2}{3} & 1 \\
  \varphi & -\frac{4}{3} & \varphi'
\end{array} \nonumber
\right)
\eea
and,
\bea
\sigma^{-1}_{3}=\left(
\begin{array}{ccc}
  \frac{2(\varphi'+2)}{5\left(\varphi-\varphi'\right)} & -\frac{4(\varphi'+2)}{5 \varphi'\left(\varphi-\varphi'\right)} & \frac{2(2\varphi'-1)}{5 \varphi' \left(\varphi-\varphi'\right)}  \\
  \frac{3}{5} & \frac{3}{5} & -\frac{3}{5}\\
  -\frac{2(\varphi+2)}{5\left(\varphi-\varphi'\right)} & \frac{4(\varphi+2)}{5 \varphi\left(\varphi-\varphi'\right)} & \frac{2(1-2\varphi)}{5 \varphi \left(\varphi-\varphi'\right)}
\end{array}
\right)=\frac{2}{5\sqrt{5}}\left(
\begin{array}{ccc}
  \varphi'+2 & -2(1-2\varphi) & (2+\varphi)  \\
  \frac{3 \sqrt{5}}{2} & \frac{3\sqrt{5}}{2} & -\frac{3\sqrt{5}}{2}\\
  -(\varphi+2) & 2 (1-2\varphi') & -(2+\varphi')
\end{array}
\right) \nonumber
\eea
Since eigenvalues of matrix $A_{3}$ are $\varphi^2,-1,\varphi'^2$, the diagonal matrix $\phi_{3}$ is,
\bea
\phi_{3} =\left(
\begin{array}{ccc}
  \varphi'^2 & 0 & 0  \\
  0 & -1 & 0\\
  0 &  0 & \varphi'^2
\end{array}
\right),
\eea
and an arbitrary $n^{th}$ power of this matrix is,
\bea
\phi^n_{3} =\left(
\begin{array}{ccc}
  (\varphi'^2)^n & 0 & 0  \\
  0 & (-1)^n & 0\\
  0 &  0 & (\varphi'^2)^n
\end{array}
\right).
\eea
Finally by using $(\ref{A3powerintermsofphiandsigma})$, $A^{n}_3=$
\begin{eqnarray}  
\frac{1}{5}\left(
  \begin{array}{ccc}
    (2F_{n}^{(2)}-3F_{n-1}^{(2)}+2(-1)^n) & (2F_{n}^{(2)}+2F_{n-1}^{(2)}+2(-1)^n) & (3F_{n}^{(2)}-2F_{n-1}^{(2)}-2(-1)^n) \\
    (F_{n}^{(2)}+F_{n-1}^{(2)}+(-1)^n) & (6F_{n}^{(2)}-4F_{n-1}^{(2)}+(-1)^n) & (4F_{n}^{(2)}-F_{n-1}^{(2)}-(-1)^n) \\
    (3F_{n}^{(2)}-2F_{n-1}^{(2)}-2(-1)^n) & (8F_{n}^{(2)}-2F_{n-1}^{(2)}-2(-1)^n) & (7F_{n}^{(2)}-3F_{n-1}^{(2)}+2(-1)^n) \\
  \end{array}
\right) \nonumber
\end{eqnarray}
is obtained.
\end{prf}

As we can expect, these results can be generalized to arbitrary matrix $A_{n+1}$. Since eigenvalues of $A_{n+1}$ are powers $\varphi^n$,$\varphi'^n$, $\ldots$, for $A^{N}_{n+1}$ eigenvalues are $\varphi^{nN}$,$\varphi'^{nN}$, \ldots But these powers can be written in terms of Fibonacci divisors as in $(\ref{phikn1})$, $(\ref{phikn2})$, and the matrix $A^{N}_{n+1}$ itself can be represented by  Fibonacci divisors $F^{(n)}_{N}$. 

For powers of matrix $A_{n+1}$ we have the following identities.

\begin{thm} \label{invarianttheorem}
Invariants of $A^{k}_{n+1}$ matrix are found as,
\bea
Tr\left( A^k_{n+1} \right)&=&\frac{F_{kn+k}}{F_{k}}=F^{(k)}_{n+1}, \label{invarianttheoremequation1} \\
{det}\left(A^k_{n+1} \right)&=&(-1)^{k \phantom{.} \frac{n(n+1)}{2}} . \label{invarianttheoremequation2}
\eea
For $k=1$, it gives
\bea
Tr\left( A_{n+1} \right)&=& F_{n+1}, \nonumber \\
{det}\left(A_{n+1} \right)&=&(-1)^{\frac{n(n+1)}{2}} . \nonumber
\eea
\end{thm}
\begin{prf}
Let's diagonalize the general matrix $A_{n+1}$ as,
\bea
\phi_{n+1}=\sigma^{-1}_{n+1}\phantom{a} A_{n+1} \phantom{a} \sigma_{n+1} \nonumber
\eea
where $\phi_{n+1}$ is diagonal and
\bea
A_{n+1}=\sigma_{n+1} \phantom{.} \phi_{n+1} \phantom{.} \sigma^{-1}_{n+1}. \nonumber
\eea
Taking the $k^{th}$ power of both sides gives,
\bea
A^{k}_{n+1}=(\sigma_{n+1} \phantom{.} \phi_{n+1} \phantom{.} \underbrace{\sigma^{-1}_{n+1})\phantom{a}(\sigma_{n+1}}_{I} \phantom{.} \phi_{n+1} \phantom{.} \sigma^{-1}_{n+1})\phantom{a}...\phantom{a}(\sigma_{n+1} \phantom{.} \phi_{n+1} \phantom{.} \underbrace{ \sigma^{-1}_{n+1})\phantom{a}(\sigma_{n+1}}_{I} \phantom{.} \phi_{n+1}  \phantom{.} \sigma^{-1}_{n+1}) \nonumber
\eea
and
\bea
A^{k}_{n+1}=\sigma_{n+1} \phantom{.} \phi^{k}_{n+1} \phantom{.} \sigma^{-1}_{n+1}. \label{An+1powerk}
\eea
By taking trace from both sides and using the cyclic permutation property of trace,
\bea
Tr(A^{k}_{n+1})=Tr\phantom{.}(\sigma_{n+1} \phantom{.} \phi^{k}_{n+1} \phantom{.} \sigma^{-1}_{n+1})=Tr\phantom{.}(\sigma^{-1}_{n+1} \phantom{.}  \sigma_{n+1} \phantom{.} \phi^{k}_{n+1})=Tr(I \phantom{a} \phi^{k}_{n+1})=Tr\phantom{.}( \phi^{k}_{n+1}) \nonumber
\eea
we get
\bea
{Tr(A^{k}_{n+1})=Tr\phantom{.}( \phi^{k}_{n+1})}. \nonumber
\eea
The eigenvalues of matrix $A_{n+1}$ in $(\ref{eigenvaluesofmatrixAn+1})$, allows one  to construct the diagonal matrix $\phi_{n+1}$ and calculate

\bea
Tr(A^{k}_{n+1})=Tr \phantom{..} \left(
                  \begin{array}{ccccccccc}
                    \varphi^n & 0 & 0 & .& . & . & 0 & 0 & 0 \\
                    0 & \varphi^{n-1} \varphi'& 0 & . & . & . & 0 & 0 & 0 \\
                    0 & 0 & \varphi^{n-2} \varphi'^2 & . & . & . & 0 & 0 & 0 \\
                    \vdots & \vdots & \vdots  & \vdots  & \vdots  & \vdots  & \vdots  & \vdots  & \vdots  \\
                    0 & 0 & 0 & . & . & . & \varphi^{2} \varphi'^{n-2} & 0 & 0 \\
                    0 & 0 & 0 & . & . & . & 0 & \varphi \varphi'^{n-1} & 0 \\
                    0 & 0 & 0 & . & . & . & 0 & 0 & \varphi'^{n} \\
                  \end{array}
                \right)^{k}. \nonumber
\eea
It gives
\bea
Tr(A^{k}_{n+1})=Tr \phantom{..}\left(
                  \begin{array}{ccccccccc}
                    (\varphi^n)^k & 0 & 0 & .  &  0 & 0 & 0 \\
                    0 & (\varphi^{n-1} \varphi')^k & 0 & .  & 0 & 0 & 0 \\
                    0 & 0  & (\varphi^{n-2} \varphi'^2)^k & . & 0 & 0 & 0  \\
                  \vdots & \vdots & \vdots  & \vdots    & \vdots  & \vdots  & \vdots  \\
                    0 & 0  & 0  & . & (\varphi^{2} \varphi'^{n-2})^k & 0 & 0 \\
                    0 & 0  & 0  &. & 0 & (\varphi \varphi'^{n-1})^k & 0 \\
                    0 & 0  & 0  &. & 0 & 0 & (\varphi'^{n})^k \\
                  \end{array}
                \right)  \nonumber \label{matrixfipowerk}
\eea
and
\be
Tr(A^{k}_{n+1})=(\varphi^n)^k+(\varphi^{n-1} \varphi')^k+\ldots +(\varphi \varphi'^{n-1})^k+(\varphi'^{n})^k, \nonumber 
\ee
or
\be
Tr(A^{k}_{n+1})=(\varphi^k)^n+(\varphi^{k})^{n-1} \varphi'^k+\ldots +\varphi^k (\varphi'^k)^{n-1}+(\varphi'^{k})^n . \nonumber
\ee
The powers $(\varphi^k)^n$ and $(\varphi'^{k})^n$ substituted from equations $(\ref{phikn1})$ and $(\ref{phikn2})$ give $Tr(A^{k}_{n+1}) =$
\bea
&=&\left(\varphi^k \phantom{.} F^{(k)}_n + (-1)^{k+1} \phantom{.} F^{(k)}_{n-1}\right)+\left(\varphi^k \phantom{.} F^{(k)}_{n-1} + (-1)^{k+1} \phantom{.} F^{(k)}_{n-2}\right)\varphi'^{k}+\ldots \nonumber \\
&&+\left(\varphi^k \phantom{.} F^{(k)}_1 + (-1)^{k+1} \phantom{.} F^{(k)}_{0}\right)(\varphi'^{k})^{n-1}+(\varphi'^{k})^{n} \nonumber  \\
&=& \varphi^k \left(F^{(k)}_n+F^{(k)}_{n-1}(\varphi'^k)+F^{(k)}_{n-2}(\varphi'^k)^2+\ldots +F^{(k)}_{1}(\varphi'^k)^{n-1}\right)  \nonumber \\ &&+(-1)^{k+1}\left(F^{(k)}_{n-1}+F^{(k)}_{n-2}(\varphi'^k)+F^{(k)}_{n-3}(\varphi'^k)^2+\ldots +F^{(k)}_{0}(\varphi'^k)^{n-1}\right)  \nonumber \\
&&+(\varphi'^k)^n \nonumber \\
&=&\varphi^k
\left(\frac{F_{kn}}{F_{k}}+\frac{F_{(n-1)k}}{F_{k}}(\varphi'^k)+\frac{F_{(n-2)k}}{F_{k}}(\varphi'^k)^2+\ldots +\frac{F_{k}}{F_{k}}(\varphi'^k)^{n-1}\right) \nonumber \\
&&+(-1)^{k+1}\left(\frac{F_{(n-1)k}}{F_{k}}+\frac{F_{(n-2)k}}{F_{k}}(\varphi'^k)+\frac{F_{(n-3)k}}{F_{k}}(\varphi'^k)^2+\ldots +\frac{F_{0}}{F_{k}}(\varphi'^k)^{n-1}\right)\nonumber\\
&&+(\varphi'^k)^n \nonumber \\
&=&\frac{F_{kn}}{F_{k}}\phantom{..} \varphi^k + \frac{F_{(n-1)k}}{F_{k}} (-1)^{k} + \frac{F_{(n-2)k}}{F_{k}} (-1)^{k} (\varphi'^k)+\ldots +\frac{F_{k}}{F_{k}}(\varphi^k) (\varphi'^k)^{n-1}\nonumber \\
&&+ \frac{F_{(n-1)k}}{F_{k}} (-1)^{k+1} + \frac{F_{(n-2)k}}{F_{k}} (-1)^{k+1}(\varphi'^k) + \frac{F_{(n-3)k}}{F_{k}} (-1)^{k+1}(\varphi'^k)^2 \nonumber \\
&&+\ldots +\frac{F_{0}}{F_{k}} (-1)^{k+1}\varphi'^{n-1} + (\varphi'^k)^n  \nonumber \\
&=&\frac{F_{kn}}{F_{k}}\phantom{..} \varphi^k + \frac{F_{(n-1)k}}{F_{k}} (-1)^{k} + \frac{F_{(n-2)k}}{F_{k}} (-1)^{k} (\varphi'^k)+\ldots  \nonumber \\
&&+ \frac{F_{(n-(n-1))k}}{F_{k}} (-1)^{k}(\varphi'^k)^{n-2} + \frac{F_{(n-1)k}}{F_{k}} (-1)^{k+1}+ \frac{F_{(n-2)k}}{F_{k}} (-1)^{k+1}(\varphi'^k) \nonumber \\
&&+ \frac{F_{(n-3)k}}{F_{k}} (-1)^{k+1}(\varphi'^k)^{2}+\ldots +\frac{F_{k}}{F_{k}} (-1)^{k+1}(\varphi'^k)^{n-2} + (\varphi'^k)^{n} \nonumber \\
&=& \frac{F_{kn}}{F_{k}}\phantom{..} \varphi^k + \frac{F_{(n-1)k}}{F_{k}} \left( (-1)^k + (-1)^{k+1} \right) \nonumber \\
&&+ \frac{F_{(n-2)k}}{F_{k}} \left( (-1)^k \varphi'^k  +(-1)^{k+1} \varphi'^k \right) + \frac{F_{(n-3)k}}{F_{k}} \left((-1)^k (\varphi'^k)^2 + (-1)^{k+1} (\varphi'^k)^2 \right)\nonumber \\
&&+\ldots +\frac{F_{k}}{F_{k}} \left( (-1)^k (\varphi'^k)^{n-2} + (-1)^{k+1} (\varphi'^k)^{n-2} \right) + (\varphi'^k)^{n} \nonumber \\
&=&\frac{F_{kn}}{F_{k}}\phantom{..} \varphi^k + \frac{F_{(n-1)k}}{F_{k}} (-1)^{k} (1+(-1)) + \frac{F_{(n-2)k}}{F_{k}} (-1)^{k} \varphi'^k (1+(-1))\nonumber \\
&&+\frac{F_{(n-3)k}}{F_{k}} (-1)^{k} (\varphi'^k)^2 (1+(-1))+\ldots +\frac{F_{k}}{F_{k}} (-1)^{k} (\varphi'^k)^{n+2} (1+(-1))\nonumber \\
&&+(\varphi'^k)^{n} \nonumber \\
&=&\frac{F_{kn}}{F_{k}} \phantom{..} \varphi^k + (\varphi'^k)^{n} \nonumber \\
&{(\ref{phikn2})}{=}&\frac{F_{kn}}{F_{k}}\phantom{..} \varphi^k + \varphi'^k \phantom{.} F^{(k)}_{n} + (-1)^{k+1} \phantom{.} F^{(k)}_{n-1} \nonumber \\
&=&\frac{F_{kn}}{F_{k}} \phantom{..} \varphi^k + \varphi'^k \phantom{.} \frac{F_{kn}}{F_{k}}  + (-1)^{k+1} \phantom{.} \frac{F_{k(n-1)}}{F_{k}} \nonumber
\eea
\bea
&=&\frac{1}{F_{k}} \left( F_{kn} \varphi^k + \varphi'^k \phantom{.} F_{kn}  + (-1)^{k+1} \phantom{.} F_{k(n-1)} \right)\nonumber \\
&=&\frac{1}{F_{k}} \frac{1}{\varphi-\varphi'} \left[ \left(\varphi^{kn}-\varphi'^{kn}\right)\varphi^{k}+\varphi'^{k} \left(\varphi^{kn}-\varphi'^{kn}\right)+ (-1)^{k+1}\left(\varphi^{(n-1)k}-\varphi'^{(n-1)k}\right) \right]   \nonumber \\
&=&\frac{1}{F_{k}} \frac{1}{\varphi-\varphi'} \left[ \varphi^{k(n+1)}-\varphi'^{kn} \varphi^{k}+\varphi'^{k} \varphi^{kn}-\varphi'^{k+kn}+ (-1)^{k+1} \varphi^{(n-1)k}-(-1)^{k+1} \varphi'^{(n-1)k} \right]   \nonumber \\
&=&\frac{1}{F_{k}} \frac{1}{\varphi-\varphi'} \bigg[ \varphi^{k(n+1)}-\varphi'^{k(n+1)}-\left(-\frac{1}{\varphi}\right)^{kn} \varphi^{k}+\left(-\frac{1}{\varphi}\right)^{k} \varphi^{kn}+(-1)^{k+1} \varphi^{(n-1)k} \nonumber \\
&&-(-1)^{k+1} \left(-\frac{1}{\varphi}\right)^{(n-1)k} \bigg]   \nonumber \\
&=&\frac{1}{F_{k}} \frac{1}{\varphi-\varphi'} \bigg[ \varphi^{k(n+1)}-\varphi'^{k(n+1)}-(-1)^{kn} \varphi^{k(1-n)}+(-1)^{k} \varphi^{k(n-1)}-(-1)^{k} \varphi^{k(n-1)} \nonumber \\
&&+(-1)^k (-1)^{k(n-1)} \varphi^{k(1-n)} \bigg]   \nonumber \\
&=&\frac{1}{F_{k}} \frac{1}{\varphi-\varphi'} \left[ \varphi^{k(n+1)}-\varphi'^{k(n+1)}-(-1)^{kn} \varphi^{k(1-n)}+(-1)^k (-1)^{kn}(-1)^{-k}  \varphi^{k(1-n)} \right]   \nonumber \\
&=&\frac{1}{F_{k}} \frac{1}{\varphi-\varphi'} \left[ \varphi^{k(n+1)}-\varphi'^{k(n+1)}-(-1)^{kn} \varphi^{k(1-n)}+(-1)^{kn} \varphi^{k(1-n)} \right]   \nonumber \\
&=&\frac{1}{F_{k}} \phantom{.} \frac{1}{\varphi-\varphi'} \left[ \varphi^{k(n+1)}-\varphi'^{k(n+1)} \right]   \nonumber \\
&=&\frac{1}{F_{k}}\phantom{.}  \frac{\varphi^{k(n+1)}-\varphi'^{k(n+1)}}{\varphi-\varphi'}   \nonumber \\
&=&\frac{1}{F_{k}} F_{k(n+1)}  \nonumber \\
&=&\frac{F_{k(n+1)}}{F_{k}} . \nonumber
\eea
To prove  the relation for ${det}\left(A^k_{n+1} \right)$, we take the determinant from both sides in $(\ref{An+1powerk})$,
\bea
{\det}\left(A^k_{n+1} \right)={\det}\left( \sigma_{n+1} \phantom{.} \phi^{k}_{n+1} \phantom{.} \sigma^{-1}_{n+1}\right).
\eea
By using property of determinants,
\bea
\det(AB)=\det(A) \det(B)
\eea
we obtain,
\bea
{\det}\left(A^k_{n+1} \right)&=&{\det}\left( \sigma_{n+1}\right) \phantom{.} {\det} \left( \phi^{k}_{n+1} \right) \phantom{.} {\det} \left( \sigma^{-1}_{n+1}   \right) \Rightarrow \nonumber
\eea
\bea
{\det}\left(A^k_{n+1} \right)&=&{\det}\left( \sigma_{n+1}\right) \phantom{.} {\det} \left( \sigma^{-1}_{n+1}   \right) \phantom{.}{\det} \left( \phi^{k}_{n+1} \right) \nonumber \Rightarrow \\
{\det}\left(A^k_{n+1} \right)&=&{\det}\left( \sigma_{n+1}  \phantom{.}\sigma^{-1}_{n+1}   \right) \phantom{.}{\det} \left( \phi^{k}_{n+1} \right) \nonumber \Rightarrow \\
{\det}\left(A^k_{n+1} \right)&=&{\det}\left(I\right) \phantom{.}{\det} \left( \phi^{k}_{n+1} \right) \nonumber  \Rightarrow\\
{\det}\left(A^k_{n+1} \right)&=&{\det} \left( \phi^{k}_{n+1} \right) .\nonumber
\eea
Since the matrix $\phi^{k}_{n+1}$ is known, the above equation becomes,
\bea
{\det}\left(A^k_{n+1} \right)&=& (\varphi^n)^k \phantom{.} (\varphi^{n-1} \varphi')^k \phantom{.} (\varphi^{n-2} \varphi'^2)^k \ldots (\varphi^{2} \varphi'^{n-2})^k \phantom{.} (\varphi \varphi'^{n-1})^k \phantom{.}(\varphi'^n)^k \nonumber \\
&=&\left(\varphi^{nk}\phantom{.}\varphi^{(n-1)k}\phantom{.}\varphi^{(n-2)k} \ldots \varphi^{2k}\phantom{.}\varphi^{k}\right)\left(\varphi'^{k}\phantom{.}\varphi'^{2k} \ldots \varphi'^{(n-2)k}\phantom{.}\varphi'^{(n-1)k}\phantom{.}\varphi'^{nk}\right) \nonumber \\
&=&\left(\varphi^{nk+(n-1)k+(n-2)k+\ldots +2k+k}\right)\left(\varphi'\phantom{.}^{k+2k+ \ldots +(n-2)k+(n-1)k+nk}\right) \nonumber \\
&=&\varphi^{k\left[n+(n-1)+(n-2)+\ldots +2+1\right]} \phantom{.} \varphi'\phantom{.}^{k\left[1+2+ \ldots +(n-2)+(n-1)+n\right]} \nonumber \\
&=& \varphi^{k\left(\frac{n(n+1)}{2}\right)} \phantom{.}  \varphi'^{k\left(\frac{n(n+1)}{2}\right)} \nonumber \\
&=& \left(\varphi^{\frac{n(n+1)}{2}}\right)^k \phantom{.} \left(\varphi'^ {\frac{n(n+1)}{2}}\right)^k  \nonumber \\
&=& \left[(\varphi \varphi')^{\frac{n(n+1)}{2}}\right]^k \nonumber \\
&{\left(\varphi \varphi'=-1\right)}{=}& (-1)^{k \phantom{.} \frac{n\phantom{.}(n+1)}{2}}.  \nonumber
\eea
\end{prf}

The above Theorem represnts Fibonacci divisors $F^{(k)}_{n+1}$ in terms of combinatorial matrix $A_{n+1}$. Quantum calculus for such divisors was constracted recently in \cite{FNK}. As was shown, it is related with several problems from hydrodynamics, quantum integrable systems and quantum information theory. This is why results of the present paper can be useful in the studies of this calculus and its applications.

\section{Acknowledgements} One of the authors (O.K.P) would like to thanks Professor Johann Cigler for attracting our attention on equivalence  of Golden binomials, introduced in \cite{golden} with Carlitz characteristic polynomials \cite{Carlitz}. This work is supported by TUBITAK grant 116F206.

\end{document}